\documentclass[]{aiaa-tc}% insert '[draft]' option to show overfull boxes
\usepackage{cite}
\usepackage{graphicx}
\usepackage{amsmath}
\usepackage{amsfonts}
\usepackage{amsthm}
\usepackage{amssymb}
\usepackage{amstext}
\usepackage{amsbsy}
\usepackage{amsopn}
\usepackage{array}
\usepackage{subfigure}
\usepackage{bm}
\usepackage{threeparttable}
\usepackage{dcolumn}
  \newcolumntype{d}{D{.}{.}{-1}}
\usepackage{hhline}
\usepackage{xspace}
\usepackage[ruled,linesnumbered]{algorithm2e}

\title{Arnoldi-based Sampling for High-dimensional Optimization using Imperfect Data}

\author{
  Jason E. Hicken\thanks{Assistant Professor, Department of Mechanical, Aerospace, and Nuclear Engineering, Member AIAA} \ and
  Anthony Ashley\thanks{Graduate Student, Department of Mechanical, Aerospace,
    and Nuclear Engineering, Student Member AIAA}    
\\
  {\normalsize\itshape
    Rensselaer Polytechnic Institute, Troy, New York, 12180}
}

 \AIAApapernumber{YEAR-NUMBER}
 \AIAAconference{Conference Name, Date, and Location}
 \AIAAcopyright{\AIAAcopyrightD{YEAR}}

 \newcommand{\ignore}[1]{}

 \newcommand{\fnc}[1]{\ensuremath{\mathcal{#1}}}

 \DeclareMathOperator{\diag}{diag}
 
 \DeclareMathOperator{\spn}{span}
 \newcommand{\etal}[0]{{\em et~al.\@}\xspace}
 \newcommand{\eg}[0]{{e.g.\@}\xspace}
 \newcommand{\ie}[0]{{i.e.\@}\xspace}

 \SetKwComment{mycomment}{}{}
 \newcommand{\algcmnt}[1]{\mycomment*[h]{\rule{3ex}{0ex}\textsf{(#1)}}}

\newcommand{\J}[0]{\fnc{J}}

\newcommand{\mat}[1]{\ensuremath{\mathsf{#1}}}
\newcommand{\gred}[0]{\ensuremath{g_{\mathsf{red}}}}
\newcommand{\xnew}[0]{\ensuremath{x_{\mathsf{new}}}}

\begin{document}

\maketitle

\begin{abstract}
  We present a sampling strategy suitable for optimization problems
  characterized by high-dimensional design spaces and noisy outputs.  Such
  outputs can arise, for example, in time-averaged objectives that depend on
  chaotic states.  The proposed sampling method is based on a generalization of
  Arnoldi's method used in Krylov iterative methods.  We show that Arnoldi-based
  sampling can effectively estimate the dominant eigenvalues of the underlying
  Hessian, even in the presence of inaccurate gradients.  This spectral
  information can be used to build a low-rank approximation of the Hessian in a
  quadratic model of the objective.  We also investigate two variants of the
  linear term in the quadratic model: one based on step averaging and one based
  on directional derivatives.  The resulting quadratic models are used in a
  trust-region optimization framework called the Stochastic Arnoldi's Method
  (SAM).  Numerical experiments highlight the potential of SAM relative to
  conventional derivative-based and derivative-free methods when the design
  space is high-dimensional and noisy.
\end{abstract}

%\section*{Nomenclature}

%\begin{tabbing}
%  XXXXX \= \kill% this line sets tab stop
%  \mb{Q}             \> conservative flow variables \\
%  $\mb{E}_{i}$       \> flux vector in direction $x_{i}$ \\
%  $\mb{\hat{E}}_{i}$ \> flux vector in direction $\xi_{i}$ \\
%  $\sigma_{1},\sigma_{2}$ \> simultaneous approximation term penalty parameters \\
%  \op{F}             \> column vector of the discrete Euler equations \\
%  \mb{q}             \> block column vector of the \mb{Q}s \\
%  \mb{f}             \> block column vector of the \mb{E}s \\
%  \mat{A}            \> discrete system Jacobian matrix \\
%  $\mat{A}_{1}$      \> first-order, modified Jacobian matrix \\
%  $A$                \> generic matrix \\
%  $S$                \> global, Schur complement matrix \\
%  $S_{i}$            \> process-local Schur complement matrix
% \end{tabbing}

%
%=======================================================================
%

\section{Introduction}\label{sec:background}

\ignore{Our premise is that gradient data, even when contaminated, provides valuable
information regarding the dominant subspaces of the objective function.  We are
not unique in this view; mention work by Paul and Qiqi here (different goal).}

Large-scale numerical simulations play a central role in contemporary aircraft
design, and, increasingly, this role goes beyond the use of simulations as
proxies for experiments.  This trend is exemplified by
differential-equation-constrained optimization (DECO), which couples simulations
with numerical optimization.  DECO has been used to optimize highly-refined
aircraft where small improvements translate into significant economic and
environmental benefits\footnote{For example, a 1\% reduction in the fuel
  consumed by the worldwide fleet would result in 7 million fewer tonnes of C02
  emitted per year~\cite{atag:2014}.}.  Moreover, DECO has the potential to
enable the design of unprecedented aircraft configurations where empirical data
is sparse and intuition is lacking.

Clearly, DECO has enormous value and potential for aerospace engineers; however,
while DECO algorithms for steady and periodic deterministic systems are
maturing, there remains a broad class of problems that cannot be optimized with
conventional algorithms.  These problems exhibit
\begin{enumerate}
\item a high-dimensional design space,
\item complex physics that must be modeled using large-scale simulations, and
\item simulation outputs, \eg lift force or total energy, that are ``imperfect''.
\end{enumerate}
In this context, \emph{imperfect} outputs are quantities of interest whose
numerical errors cannot be eliminated, at least in practice; consequently, such
outputs fail to meet the underlying assumptions of conventional gradient-based
optimization algorithms.  The word imperfect is intentionally chosen to
distinguish these errors from more traditional numerical errors that can,
usually, be estimated and effectively reduced.  In the following sections, we
briefly elaborate on two such sources of imperfect data.

\subsection{Time-averaged Outputs from Chaotic Systems}

The outputs of interest in many engineering systems are time-averages of chaotic
solutions.  Relevant examples include the lift and drag on aerodynamic
bodies~\cite{spalart:2009}, the energy produced by a fusion
reactor~\cite{debaar:2005}, and the (phase-averaged) pressure in an
internal-combustion engine~\cite{enotiadis:1990}.

\begin{figure*}[tbp]
  \begin{center}
    \subfigure[\label{fig:lorenz_trajectory}]{\includegraphics[width=0.49\textwidth]{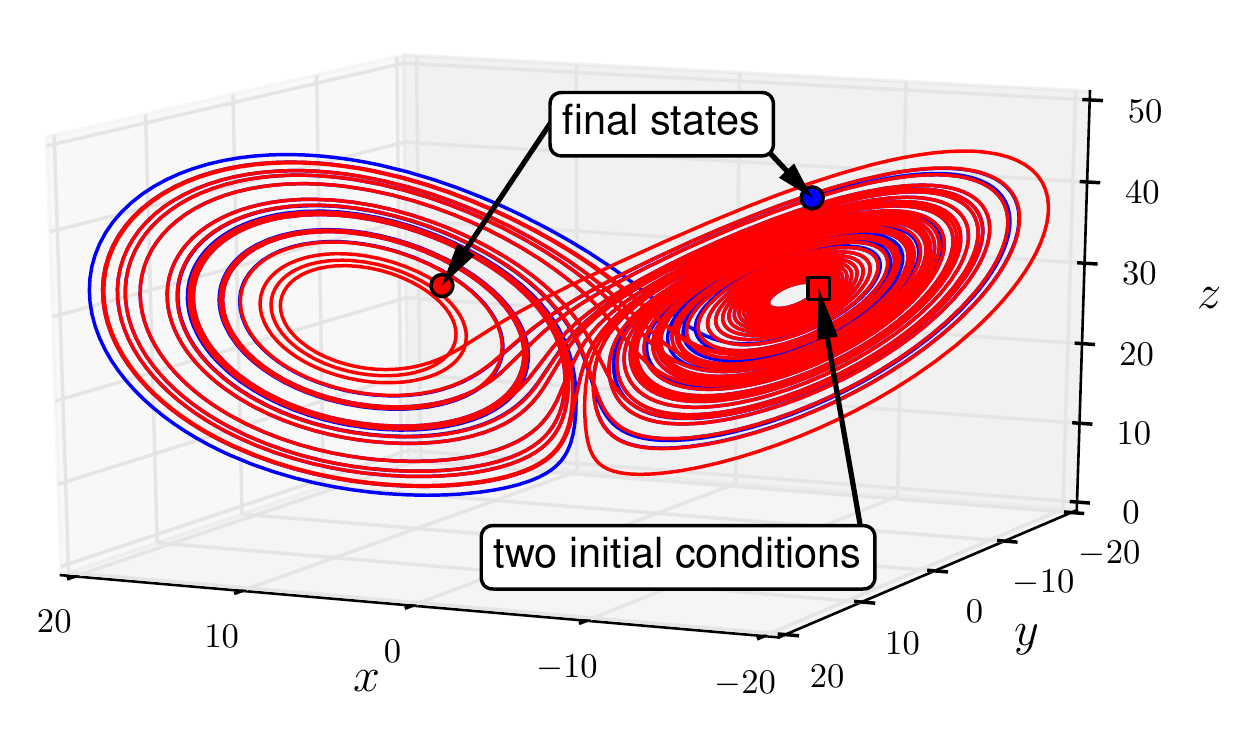}}
    \subfigure[\label{fig:lorenz_obj}]{%
      \includegraphics[width=0.49\textwidth]{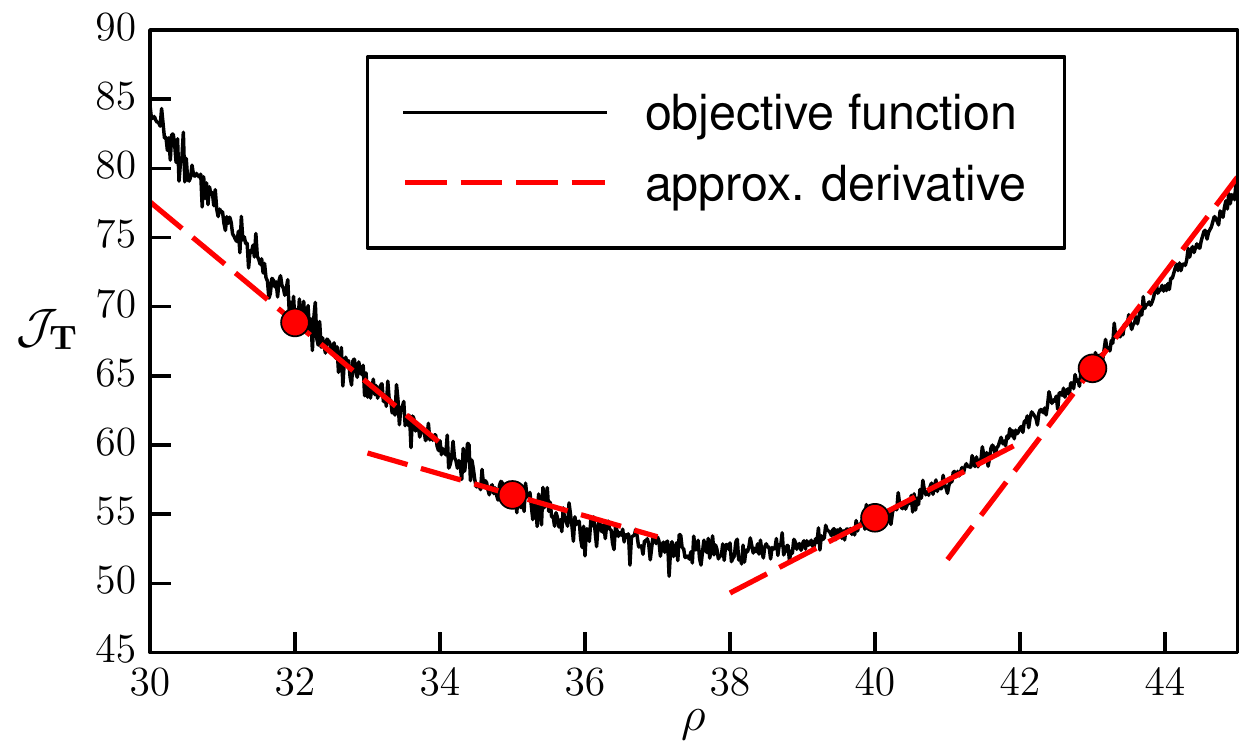}}
    \caption[]{Example trajectories of the Lorenz DE (left) illustrating
      sensitive dependence on initial conditions.  A time-averaged objective
      (right) exhibits high-frequency fluctuations; moreover, only approximate
      derivatives can by computed, in this case using an ensemble adjoint.}
  \end{center}
\end{figure*}

Chaotic systems are characterized by a sensitive dependence on initial
conditions; the distance between two states that are ``close'' initially will
exponentially increase with time.  This is illustrated in
Figure~\ref{fig:lorenz_trajectory} using the Lorenz DE~\cite{lorenz:1963} and
two initial conditions that satisfy $\| \Delta \bm{x}_{0} \| \leq 10^{-10}$.

In addition to the initial conditions, chaotic systems are sensitive to other
parameters.  This has significant implications for time-averaged outputs, which
we illustrate using the Lorenz system and the objective function
\begin{equation*}
  \J_{T}(\rho) = \frac{1}{2T} \int_{0}^{T} \left( z(t,\rho) - z_{\text{targ}}\right)^{2} \, dt,
\end{equation*}
where $z(t,\rho)$ is one of the Lorenz state variables, $z_{\text{targ}} = 35$, and
$T > 0$ is the period of integration.  The design variable here is $\rho$, a
parameter in the Lorenz DE.\ignore{  The Lorenz-DE parameters $\sigma$ and $\beta$ are
fixed at 10 and $8/3$, respectively.}

\ignore{
Time-averaged objective functions, similar to $\J_{T}(\rho)$, frequently appear
in engineering analyses.  For example, the time-averaged lift/drag on an
aircraft or the time-averaged energy produced by a reactor.  For such objectives
we would ideally like $T \rightarrow \infty$; in practice we must truncate the
simulation with a finite $T$.}

Figure~\ref{fig:lorenz_obj} plots the Lorenz objective $\J_{T}$ versus the
parameter $\rho$ for an averaging period of $T = 400$.  High-frequency
oscillations can be observed in the objective function $\J_{T}$.  In theory, we
could eliminate these fluctuations by integrating over an infinite time horizon,
but, in practice, we must truncate the simulation at a finite $T$.  The
oscillations in $\J_{T}$ reflect the sensitivity of the Lorenz DE to changes in
$\rho$, and they hint at the difficulties of using gradient-based optimization.
Indeed, as $T \rightarrow \infty$ the gradient of $\J_{T}$
diverges~\cite{lea:2000} despite the fact that the objective itself converges.

%Note on derivative-free?

Researchers have proposed methods to compute derivatives of objectives that
depend on chaotic systems, such as the ensemble adjoint~\cite{lea:2000,
  eyink:2004, ashley:2014} and least-squares adjoint~\cite{wang:2013b}.  These
methods share one shortcoming: \emph{they produce estimates of the derivatives
  only}.  Figure~\ref{fig:lorenz_obj} illustrates some derivatives estimated
using the ensemble adjoint.  Although the derivatives capture the general slope
of the objective, errors are clearly visible.  Such errors are incompatible with
conventional gradient-based optimization.

\subsection{Uncertainty Propagation in High-dimensional Input Spaces}

Numerical simulations contain many sources of uncertainty.  For example,
turbulence and combustion models introduce errors with respect to the true
physics, the operating conditions of a system are not deterministic (\eg the
Reynolds number is not known precisely), and perturbations to the intended
design are introduced by the manufacturing process.  The field of optimization
under uncertainty (OUU), a branch of uncertainty quantification (UQ), seeks to
account for such uncertainties during the optimization process.  OUU can help
ensure good performance over a range of parameters, and it can help reduce the
probability of failure.

Typically, at each iteration of an OUU, uncertainties in the input parameters
must be propagated through the numerical simulation to deduce uncertainties in
the outputs, \eg the standard deviation in a force or energy.  Uncertainty
propagation is especially challenging when there are a large ($> 100$) number of
uncertain input parameters.  These situations require approximations to cope
with the \emph{``curse of dimensionality''}.

Several parameter-selection and model-reduction strategies have been proposed
for the propagation of stochastic, or aleatory, uncertainties in
high-dimensional input spaces~\cite{veroy:2005, buithanh:2008, bashir:2008}.
Methods of propagating epistemic uncertainties, \ie model errors, in
high-dimensional spaces have also been investigated~\cite{lockwood:2012,
  boopathy:2014}.  If these methods are to be used for high-dimensional OUU
problems, derivative information will be needed; however,
\emph{differentiating these propagation methods accurately is not practical}.
Therefore, an optimization algorithm must account for the discrepancy between
the outputs and their derivatives.

\ignore{
\begin{itemize}
\item Several parameter-selection and model-reduction strategies have been
  proposed for the propagation of stochastic, or aleatory, uncertainties in
  high-dimensional input spaces; for example,~\cite{veroy:2005, buithanh:2008,
    bashir:2008}.  These methods have shown promise for linear state equations.
  Parameter selection in the nonlinear case is significantly more
  challenging~\cite{smith:2014}, and errors between the reduced model and the
  ``true'' model are likely to be significantly higher for the same
  computationally cost.
\item For propagation of epistemic uncertainties, Lockwood
  \etal~\cite{lockwood:2012} propose using a bound-constrained optimization to
  find the worst outcome for an output over the specified domain of the
  epistemic variables.  To avoid a nested optimization in OUU, Boopathy and
  Rumpfkeil~\cite{boopathy:2014} used the derivatives at the midpoint of the
  epistemic domain as a proxy for the derivative at the min/max.  This
  introduces an error in the gradient needed for OUU.
\end{itemize}
}

\begin{figure*}[tbp]
  \begin{center}
    \subfigure[\label{fig:rosenbrock}]{\includegraphics[height=0.2\textheight]{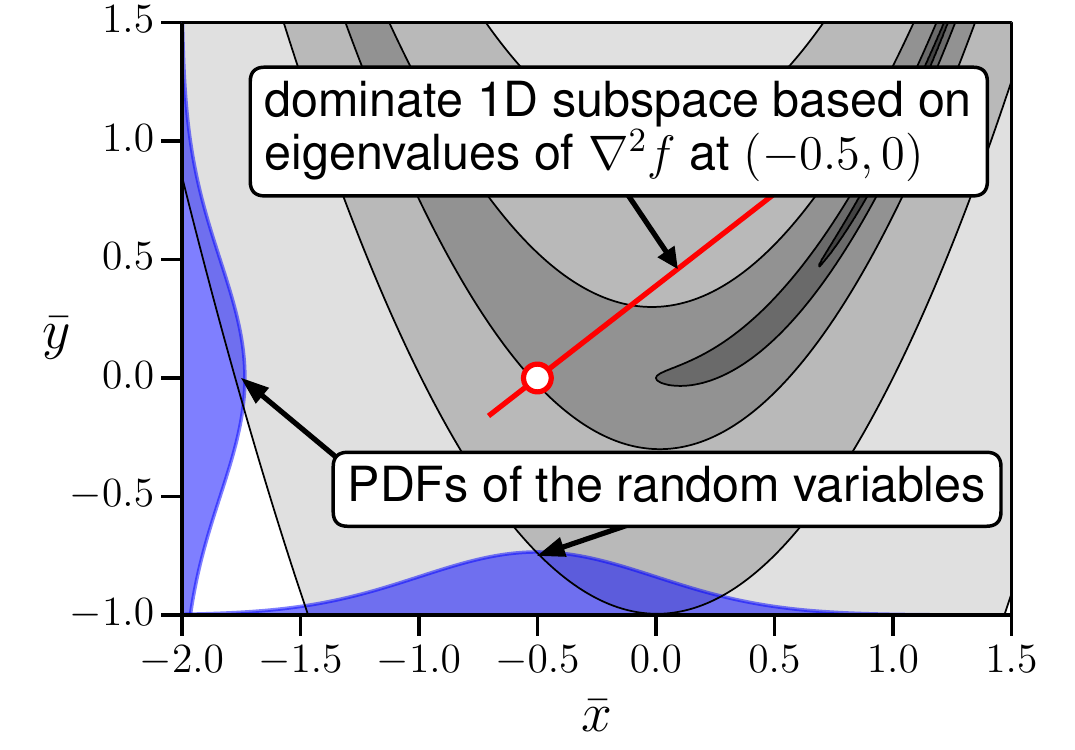}}
    \subfigure[\label{fig:UQ_subspace}]{\includegraphics[height=0.2\textheight]{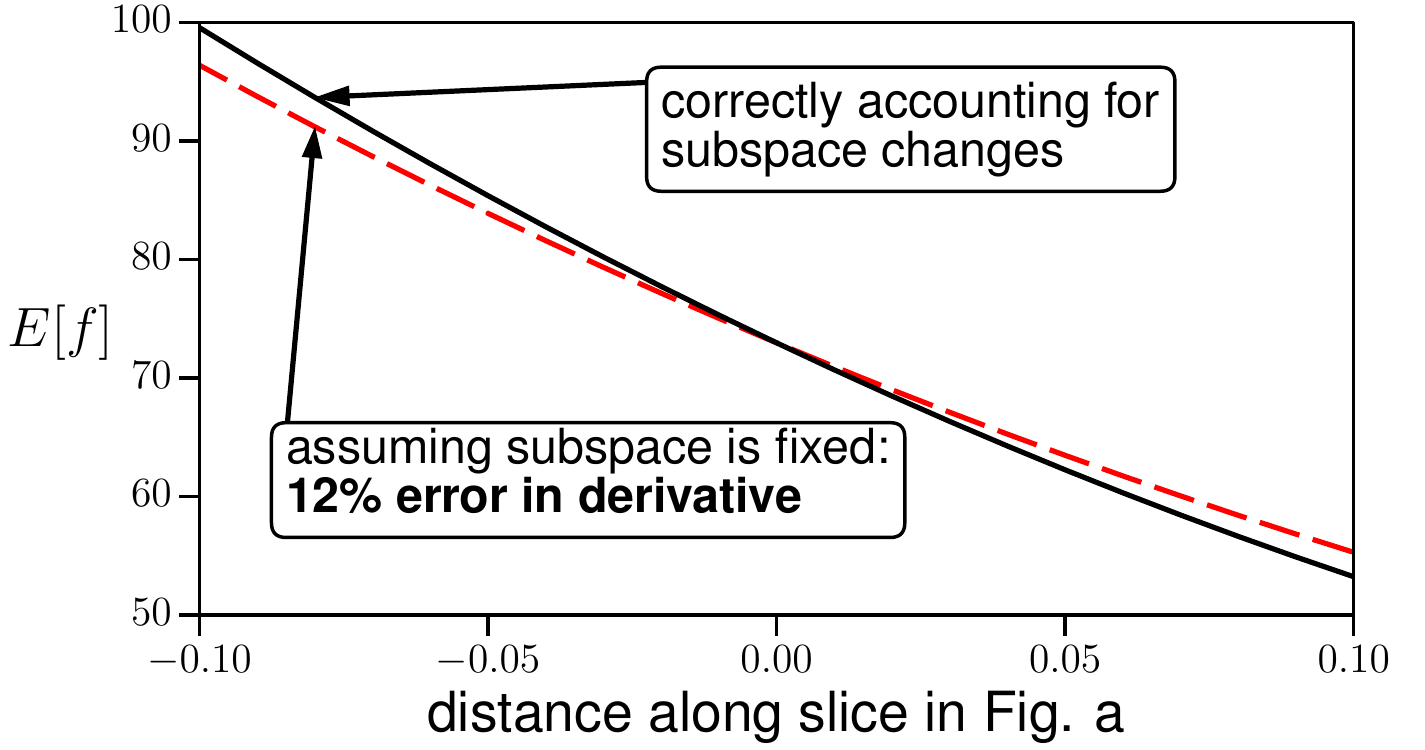}}
    \caption[]{Uncertainty propagation can be made tractable for large numbers
      of (aleatory) inputs by identifying dominant subspaces using, \eg, the
      Hessian (left).  However, accurately differentiating such
      dimension-reduction strategies is usually not practical, and errors are
      inevitable (right).}
  \end{center}
\end{figure*}

To illustrate why dimension-reduction strategies are difficult to differentiate,
consider computing the expected value, $E[f(x)]$, of the Rosenbrock function
\begin{equation*}
  f(\bm{x}+\bm{\xi}) = \left[1 - (x+\xi)\right]^{2} + 100\left[(y-\eta) - (x+\xi)^{2}\right]^{2},
\end{equation*}
where $\bm{x} = (x,y)^{T} \in \mathbb{R}^{2}$ are the design variables and $\xi,
\eta \in \mathcal{N}(0,0.5^2)$ are Gaussian random variables.  In this example,
we use the dominant eigenvalue of the Hessian to define a subspace for dimension
reduction.  Figure~\ref{fig:rosenbrock} shows one such subspace defined at
$\bm{x} = (-0.5, 0)^{T}$.

Figure~\ref{fig:UQ_subspace} plots the expected value obtained by (correctly)
accounting for changes in the subspace.  The figure also plots the expected
value obtained by fixing the subspace, which reflects how derivatives would be
computed for most dimension-reduction strategies.  Fixing the subspace results
in a 12\% error in the derivative.  Note that, to correctly compute the
derivative in the present example, we would need to differentiate eigenvectors
of the Hessian, which would be impractical for high-dimensional outputs based on
simulations.

\ignore{
Moreover, while the reduced-space objective is smooth in
Figure~\ref{fig:UQ_subspace}, this may not be the case in general. In
particular, when (random) importance sampling is used~\cite{ghattas}, the
estimate of $E[f]$ and its derivative will have fluctuations due to the
randomness.
}

\subsection{Other Sources of ``Imperfect'' Outputs}

Time-averaging in chaotic systems and uncertainty propagation in
high-dimensional spaces are two important examples that produce imperfect data
and motivate the current work.  However, inconsistencies between the outputs and
their derivatives can also arise in other applications of DECO.
\begin{enumerate}
\item Cut-cell~\cite{aftosmis:1997, aftosmis:1998, fidkowski:2007} and
  immersed-boundary methods~\cite{peskin:1972, mittal:2005, griffith:2005} are
  popular and effective methods for numerically solving DEs on complex, moving
  geometries.  In geometry optimization, these methods produce discontinuities
  in the design space as the mesh-topology and/or stencil is updated.  A similar
  problem exists when unstructured grids are regenerated during shape
  optimization.
\item The differentiate-then-discretize~\cite{gunzburger:2003}, or
  ``continuous''~\cite{giles:2000}, adjoint yields derivatives that are
  inconsistent with the discretized output.  See, for
  example,~\cite{collis:2002, hicken:dual2014}.
\item Incomplete sensitivities have been proposed~\cite{mohammadi:2004,
  derakhshan:2008} to reduce the computational cost of computing the gradient.
  The terms that are dropped or estimated in these approaches necessarily result
  in approximate derivatives.
\item A multi-fidelity approach may be used in which a higher fidelity tool
  produces the output but a lower fidelity method is used to estimate gradients
  \cite{balabanov:2004}.
\end{enumerate}
The inaccuracies 1 and 2 described above are discretization errors that can be
reduced through mesh refinement.  However, the mesh refinement needed to
sufficiently reduce these errors may not be possible in practice, and in these
cases the algorithms proposed below could be helpful.

\subsection{Imperfect Data and the State-of-the-art}\label{sec:limitations}

The examples above highlight important emerging applications ---
high-dimensional uncertainty propagation and simulation of chaotic systems ---
in which outputs of interest and their derivatives exhibit some form of
inaccuracy.  If we want to optimize these outputs, can we rely on
state-of-the-art optimization methods?

For low-dimensional design spaces, derivative-free optimization methods could be
effective for the applications described above.  These methods include
deterministic interpolation/regression-based surrogate models~\cite{conn:1998,
  powell:2002, conn:2009}, the Nelder-Mead simplex method~\cite{nelder:1965},
and genetic algorithms~\cite{holland:1975, hajela:1990}.  However, ``it is
usually not reasonable to try and optimize problems with more than a few tens of
variables...'' with derivative-free methods~\cite{conn:2009}.

Derivative-based algorithms are highly scalable, making them ideally suited for
optimization of smooth functions in large design spaces.  Although the challenge
of differentiating the simulation software is a potential drawback,
algorithmic-differentiation~\cite{griewank:2008} has made this task easier.
Unfortunately, derivative-based algorithms are also not suitable for the target
applications, because they require sufficient accuracy in the data and
consistency between the outputs and their derivatives.

Stochastic approximation (SA) algorithms~\cite{robbins:1951, kiefer:1952,
  spall:1992, spall:2005} can be used to optimize functions whose evaluation
contains noise.  Some forms of SA also permit the use of noisy
gradients~\cite{robbins:1951, spall:2009}.  There are both theoretical and
practical issues with applying SA to the applications described earlier.
Theoretically, the ``noise'' must be unbiased and independent of the design
space~\cite{spall:2005}, and these assumptions are not met by the target
applications.  From a practical perspective, SA algorithms tend to use many
``low quality'' iterations to ensure (probabilistic) convergence, so function
and gradient evaluations must be inexpensive.  This requirement is also not
fulfilled by the DE-constrained optimization problems under consideration.

Reduced-order models (ROM) offer a distinct approach from derivative-free and
derivative-based optimization.  Rather than tackling the DE-based optimization
directly, ROM methods first seek a simplified, and presumably less expensive,
model for the DE simulation.  Subsequently, the optimization is performed using
the ROM as a surrogate for the full DE model.  Methods in this class include
projection-based approaches like proper-orthogonal
decomposition~\cite{pearson:1901, willcox:2002, rowley:2004, volkwein:2011,
  amsallem:2012}.  These approaches typically fix the design parameters while
constructing the ROM, although a method that accounts for parameter dependence
was recently proposed by Lieberman \etal~\cite{lieberman:2010} in the case of
steady, linear DEs.  However, building parametric ROMs for nonlinear and/or
unsteady DEs remains an active and challenging area of research.

\section{Arnoldi Sampling}

The primary contribution of this work is a sampling procedure that is suitable
for high-dimensional input spaces when the gradient is available, but
potentially inaccurate.  The sampling procedure is based on Arnoldi's method,
which we briefly review below.

\subsection{Arnoldi's Method}

Arnoldi's method is found in Krylov subspace methods for spectral analysis and
solving linear systems; see, for example, \cite{saad:1993} and \cite{saad:2003}
and the references therein.  In those applications, Arnoldi's method is used to
construct an orthogonal basis for the Krylov subspace
$\mathcal{K}_{m}(\mat{A},z) \equiv \spn\{z, \mat{A}z, \mat{A}^{2}z, \cdots,
\mat{A}^{m-1}z\}$, where $\mat{A}$ is a square matrix and $z$ is a vector.

A version of Arnoldi's method based on modified Gram-Schmidt is provided in
Algorithm~\ref{alg:Arnoldi} for reference.  An important feature of Arnoldi's
method is that $\mat{A}$ is not required explicitly: only matrix-vector products
of the form $\mat{A} z_{j}$ are needed.  We exploit this aspect of the algorithm
in the proposed sampling procedure.

\begin{algorithm}[htbp]\DontPrintSemicolon
  \KwData{$z_{1} \in \mathbb{R}^{n}$ such that $\|z_{1}\| = 1$}
  \KwResult{$\mat{Z} = \left[ z_{1}, z_{2}, \ldots, z_{m+1}\right]$, 
    an orthonormal basis for $\mathcal{K}_{m+1}(\mat{A},z_{1})$}
  \BlankLine
  \For{$j = 1,2,\ldots,m$}{
    $z_{j+1} = \mat{A} z_{j}$\;
    \For(\algcmnt{Modified Gram-Schmidt}){$i=1,\ldots,j$}{
      $h_{i,j} = z_{j+1}^{T} z_{j}$\;
      $z_{j+1} \leftarrow z_{j+1} - h_{i,j} z_{j}$\;
    }
    compute $h_{j+1,j} = \| z_{j+1}\|$\;
    \lIf(\algcmnt{check for breakdown}){$\|h_{j+1,j}\| = 0$}{return}
    $z_{j+1} \leftarrow z_{j+1}/h_{j+1,j}$\;    
  }
  \caption{Arnoldi's method.\label{alg:Arnoldi}}
\end{algorithm}

To understand the origins of the proposed sampling procedure, it is helpful to
review how Arnoldi's method is traditionally used in the context of
optimization.  Efficient optimization algorithms for $C^{2}$ objectives apply
Newton's method, or a quasi-Newton method, to the first-order optimality
conditions.  For unconstrained convex problems, Newton's method produces systems
of the form
\begin{equation*}
  \mat{W} p = - g,
\end{equation*}
where $p$ is a trial step, $g = \nabla f$ is the gradient of the objective, and
$\mat{W} = \nabla^{2} f$ is the Hessian, or an approximation to it.  When
Algorithm~\ref{alg:Arnoldi} is used in this setting, $\mat{A}$ becomes $\mat{W}$
and Arnoldi's method reduces to the symmetric Lanczos algorithm that appears in
the conjugate gradient method.

Thus, when applied to optimization problems, Arnoldi's method requires
Hessian-vector products.  These products can be computed in several ways,
including algorithmic differentiation~\cite{griewank:2008} and, when PDEs are
involved, second-order adjoints~\cite{wang:1992, borzi:2011,
  hicken:inexact2014}.  The products can also be computed using a
finite-difference approximation applied to the gradient, since
\begin{equation*}
  \left[\nabla^{2} f\right] z_{j} = \lim_{\epsilon \rightarrow 0} \frac{g(x +
    \epsilon z_{j}) - g(x)}{\epsilon}.
\end{equation*}
\ignore{In practice, a finite $\epsilon$ must be chosen, requiring a careful trade-off
between truncation and round-off errors.  For this reason, the finite-difference
approximation is not the favored method when other options are available.}

\subsection{Arnoldi Sampling}

Arnoldi's method can be transformed into a sampling procedure by recognizing
that the Hessian-vector products selected by the algorithm represent
infinitesimal samples.  If these infinitesimal perturbations are made finite,
they can be used as a sampling procedure.  In other words, the sample locations
are defined by
\begin{equation*}
  x_{j} = x_{0} + \alpha z_{j},
\end{equation*}
where $\alpha > 0$ is the sample radius, and the $z_{j}$ are defined by
Arnoldi's method with $\mat{A} z_{j}$ replaced with\linebreak
$\left[g(x_{0}+\alpha z_{j}) - g(x_{0})\right]/\alpha$.  This Arnoldi sampling
procedure is listed in Algorithm~\ref{alg:Arnoldi_sample}.

In addition to providing the sample locations and sampled data,
Algorithm~\ref{alg:Arnoldi_sample} also produces approximations to the
eigenvalues and eigenvectors of the Hessian; see lines~\ref{alg:eigen_decomp} and
\ref{alg:compute_V}.  This approximation is based on iterative eigenvalue
methods~\cite{arnoldi:1951,saad:1993}.

The eigen-decomposition in Arnoldi sampling uses the symmetric part of
$\mat{H}_{m}$, the $m\times m$ upper Hessenberg matrix composed of the
$h_{i,j}$.  In contrast, the conventional Arnoldi's method for spectral analysis
uses the matrix $\mat{H}_{m}$ itself.  We use the symmetric part, because
$\mat{H}_{m}$ reduces to a symmetric matrix when the gradients are accurate and
$\alpha \rightarrow 0$.  This also avoids the generation of a complex spectrum,
which would not be appropriate in the context of a Hessian.

\ignore{ 

Why does this work? --> Can approximate the dominant eigenmodes of the
Hessian, even in the presence of errors in the gradient.

}

When the gradients are accurate and $\alpha$ is chosen suitably, Arnoldi
sampling reduces to a Arnoldi's method with finite-difference approximations for
the Hessian-vector products.  Like finite-difference approximations, $\alpha$
must be chosen carefully to achieve optimal performance.  Unlike
finite-difference approximations, our proposed method benefits from a relatively
large sampling radius, one that would normally cause undesirable truncation
errors in a finite-difference approximation.

\begin{algorithm}[htbp]\DontPrintSemicolon
  \KwData{$m > 0$, $\alpha > 0$, $x_{0}$, $f_{0} = f(x_{0})$ and $g_{0} = g(x_{0})$}
  \KwResult{sampled locations $\mat{X}_{m+1} = \left[x_{0}, x_{1}, x_{2},
      \ldots,x_{m}\right]$\;
    \hspace{4em}sampled function values $\mat{F}_{m+1} = \left[f_{0},
      f_{1}, \ldots, f_{m}\right]$\;
    \hspace{4em}sampled gradient values $\mat{G}_{m+1} =
    \left[g_{0}, g_{1}, \ldots, g_{m}\right]$\;
    \hspace{4em}approximate eigenvalues,
    $\Lambda_{m} = \left[\lambda_{1}, \lambda_{2}, \ldots, \lambda_{m} \right]$, and\;
    \hspace{4em}approximate eigenvectors of
    the Hessian $\mat{V}_{m} = \left[ v_{1}, v_{2}, \ldots, v_{m} \right]$}
  \BlankLine
  set $z_{1} = -g_{0}/\|g_{0}\|$\;
  \For{$j = 1,2,\ldots,m$}{
    set $x_{j} = x_{0} + \alpha z_{j}$\;
    sample $f_{j} = f(x_{j})$ and $g_{j} = g(x_{j})$\;
    compute $z_{j+1} = (g_{j} - g_{0})/\alpha$\;
    \For(\algcmnt{Modified Gram-Schmidt}){$i=1,\ldots,j$}{
      $h_{i,j} = z_{j+1}^{T} z_{j}$\;
      $z_{j+1} \leftarrow z_{j+1} - h_{i,j} z_{j}$\;
    }
    compute $h_{j+1,j} = \| z_{j+1}\|$\;
    \If(\algcmnt{check for breakdown}){$\|h_{j+1,j}\| = 0$}{%
      set $m = j$ and break\;
    }
    $z_{j+1} \leftarrow z_{j+1}/h_{j+1,j}$\;
  }
  compute eigen-decomposition of the symmetric part of $\mat{H}_{m}$, \ie $\displaystyle
    \frac{1}{2}\left[\mat{H}_{m} + \mat{H}_{m}^{T}\right] \tilde{\mat{V}}_{m} =
    \tilde{\mat{V}}_{m} \Lambda_{m}.$ \label{alg:eigen_decomp}\;
  compute the approximate eigenvectors $\mat{V}_{m} = \mat{Z}_{m} \tilde{\mat{V}}_{m}$
  \label{alg:compute_V}\;
  \caption{Arnoldi sampling.\label{alg:Arnoldi_sample}}
\end{algorithm}

\subsection{Eigenvalue Accuracy Study}

We use the dominant eigenvalues produced by the Arnoldi sampling procedure to
build a quadratic model for optimization; therefore, it is important to verify
and quantify the accuracy of the approximate eigenvalues.  To this end, we apply
Algorithm~\ref{alg:Arnoldi_sample} to a set of quadratic objectives with
different eigenvalue distributions.  In particular, the objective is defined by
\begin{equation}
  F(x) = x^{T} \mat{E} \Sigma \mat{E}^{T} x, \label{eq:quad_obj}
\end{equation}
where $\mat{E}$ denotes the $2^{p}\times 2^{p}$ orthonormalized Hadamard matrix,
whose columns are the synthetic eigenvectors.  The diagonal matrix $\Sigma \in 
\mathbb{R}^{2^p \times 2^p}$ holds the synthetic eigenvalues, given by
\begin{equation*}
  \Sigma_{i,i} = \frac{1}{i^q},
\end{equation*}
where $q = \frac{1}{2}, 1$, or $2$.

For this investigation, we consider $n = 256$ variables.  We run Arnoldi
sampling with $m=16$ iterations and a sample radius of $\alpha = 1$.  The
initial point about which samples are taken is $(x_{0})_{i} = \sin(i), i =
1,2,\ldots,256$.  Noise with a Gaussian distribution is added to each component
of the gradient.  The noise has mean zero and its standard deviation is either
0.5\%, 2.5\% or 5\%, relative to the norm of the initial gradient at $x_{0}$.

\begin{figure}[tbp]
  \begin{center}
  \subfigure[0.5\% noise \label{fig:eigvalues_noise01}]{%
    \includegraphics[width=0.9\textwidth]{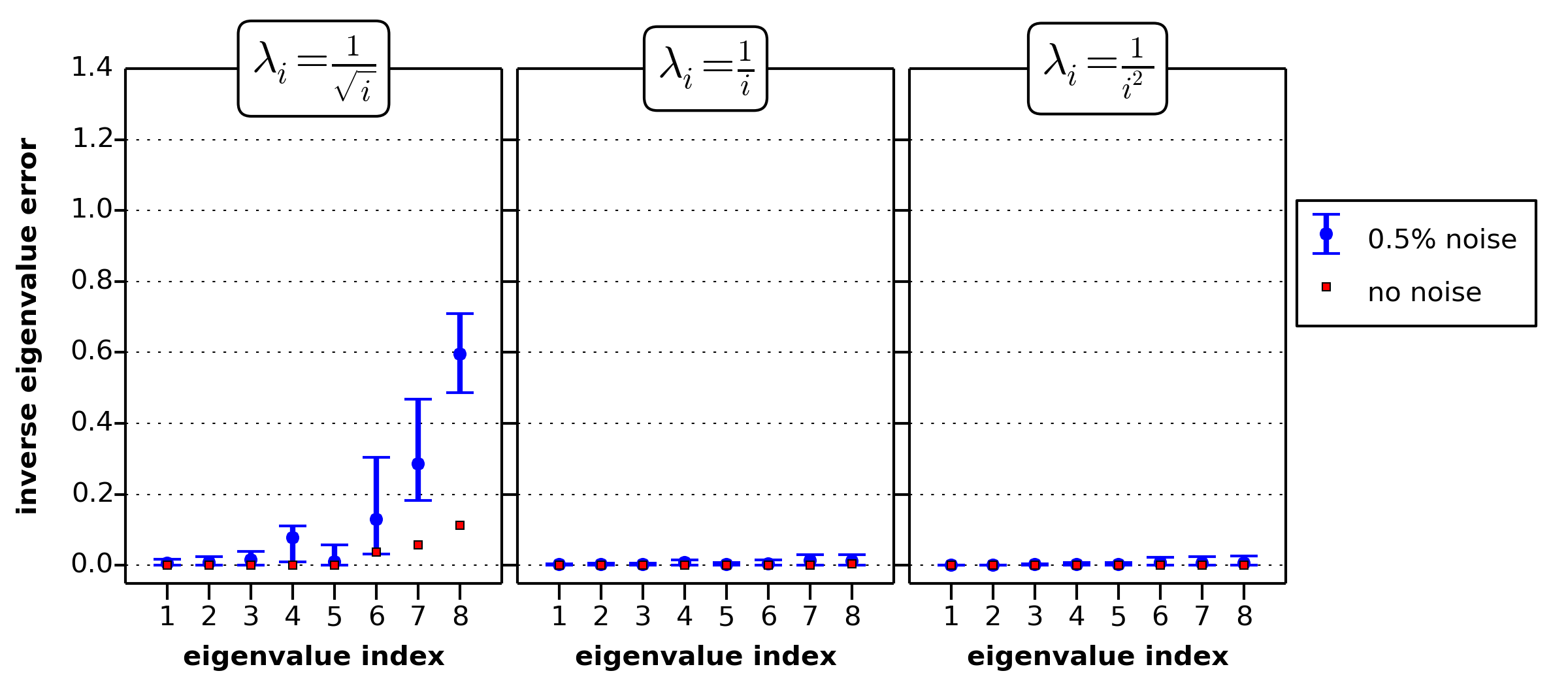}}\\
  \subfigure[2.5\% noise \label{fig:eigvalues_noise05}]{%
        \includegraphics[width=0.9\textwidth]{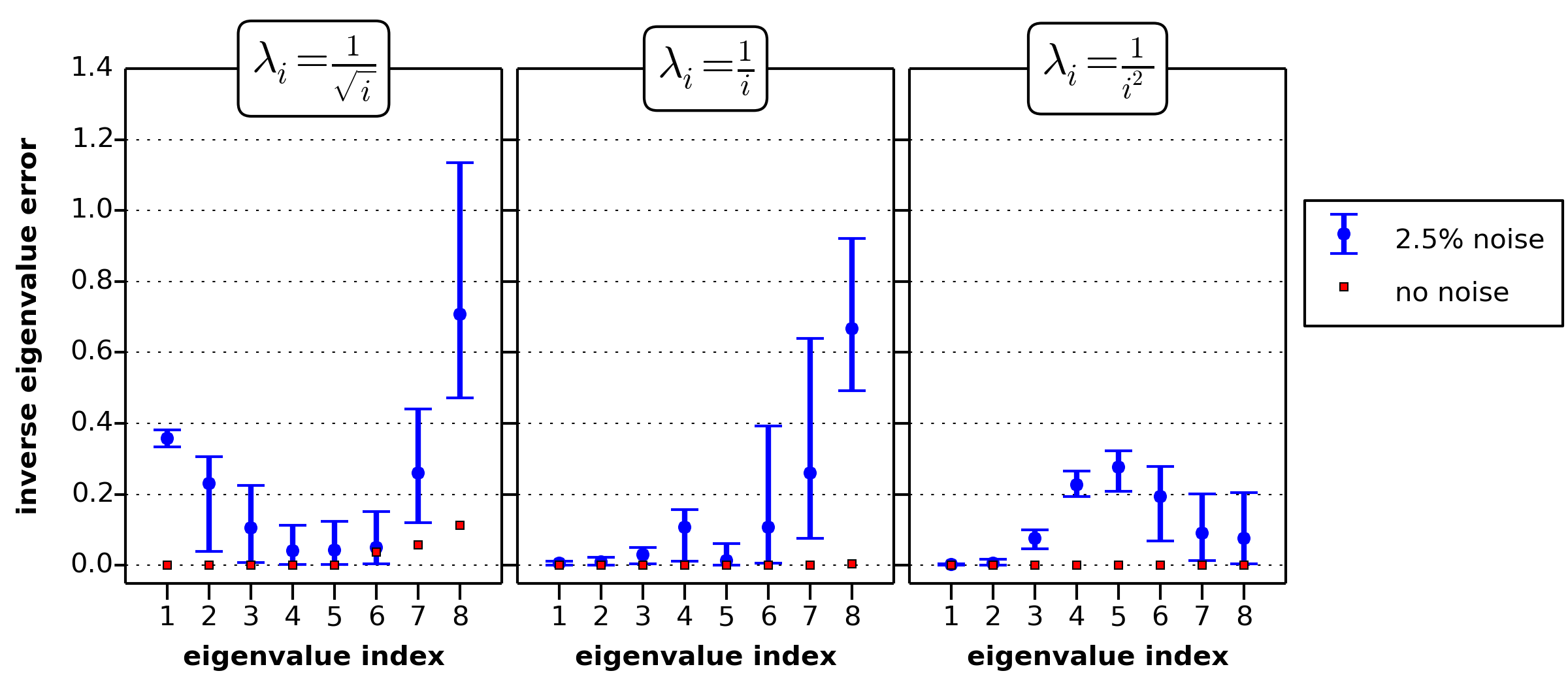}}\\
  \subfigure[5\% noise \label{fig:eigvalues_noize10}]{%
        \includegraphics[width=0.9\textwidth]{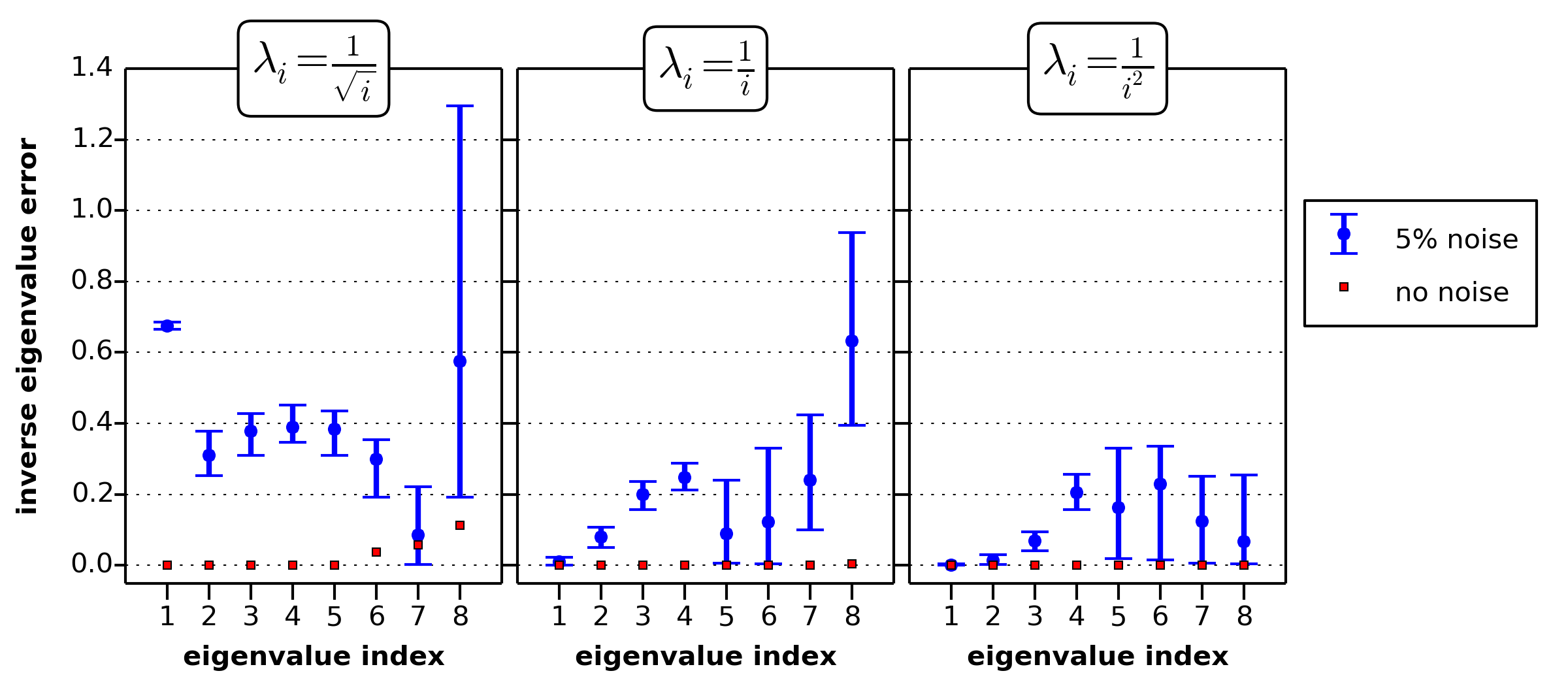}}
  \end{center}
  \caption{Relative accuracy in the inverse of the eigenvalues approximated
    using Arnoldi sampling with varying noise
    magnitude.  \label{fig:eigvalues}}
\end{figure}

Figure~\ref{fig:eigvalues} shows the relative error in the first 8 estimated
inverse eigenvalues:
\begin{equation*}
  \textsf{error}_{i} = \left| \tilde{\lambda}_{i}^{-1} - \lambda_{i}^{-1} \right| / \lambda_{i}^{-1} = \left| \lambda_{i} / \tilde{\lambda}_{i} - 1 \right|,
\end{equation*}
where $\tilde{\lambda}_{i}$ is the $i^{\text{th}}$ estimated eigenvalue and
$\lambda_{i}$ is the exact eigenvalue.  Recall that the eigenvalue distributions
are given by $\{1/\sqrt{i}\}_{i=1}^{256}$, $\{1/i\}_{i=1}^{256}$, and
$\{1/i^2\}_{i=1}^{256}$; the results for the three distributions are plotted
left to right in the figure.  The blue dot in the figures denotes the median
value of the error over 100 runs, and the lower and upper bars represent the
0.025 and 0.975 probability quantiles, respectively.  The magnitude of the noise
increases from 0.5\% in the top set of figures to 5\% in the bottom set of
figures.  For reference, the red squares denote the error in the estimated
eigenvalues when no noise is present.

Overall, the plots suggest that the dominant eigenvalues estimated by Arnoldi
sampling are resilient over a range of noise magnitudes.  In addition, the
results indicate that spectra whose energy decays quickly are less affected by
noise; nevertheless, it is important to point out that relative errors in
eigenvalues with small magnitude can have a pronounced impact on the
optimization steps.

The results in Figure~\ref{fig:eigvalues} also suggest that Arnoldi sampling
should not be used with relative noise above 5\%; however, this conclusion
pertains to the chosen radius of $\alpha = 1$ only, since an increase in
$\alpha$ tends to reduce the impact of noise.  On the other hand, an increase in
$\alpha$ risks an increase in truncation error for strongly nonlinear functions.

\section{Stochastic Arnoldi's Method}

In this section, we describe how Arnoldi sampling can be used in an optimization
framework.  The key idea is to form a quadratic model based on the estimated
eigenvalues and eigenvectors, and seek a candidate step $p$ in the span of the
eigenvectors:
\begin{equation}\label{eq:model}
    \min_{p \in \spn(\mat{V})} q(p) \equiv \bar{f} + \bar{g}^{T} p + \frac{1}{2} p^{T} \mat{V}
    \Lambda \mat{V}^{T} p, \qquad
    \textsf{subject to}\quad \| p \| \leq \Delta,
\end{equation}
where the diagonal matrix $\Lambda \in \mathbb{R}^{r\times r}$ holds the $r < m$
largest (in magnitude) eigenvalue estimates from Arnoldi sampling, and $\mat{V}
\in \mathbb{R}^{n\times r}$ are the corresponding eigenvector estimates.  Note
that \eqref{eq:model} includes a trust-radius constraint, with radius $\Delta >
0$, to handle the possibility of nonconvex models.

We have not yet specified how $\bar{f}$ and $\bar{g}$ are estimated in the
trust-region subproblem~\eqref{eq:model}.  While the function value has no
influence on the solution of $p$, it does play a role in the step-acceptance
procedure described later.  The linear term $\bar{g}$, on the other hand, does
impact $p$.

\begin{figure}[tbp]
  \begin{center}
  \subfigure[accurate gradient]{%
    \includegraphics[width=0.49\textwidth]{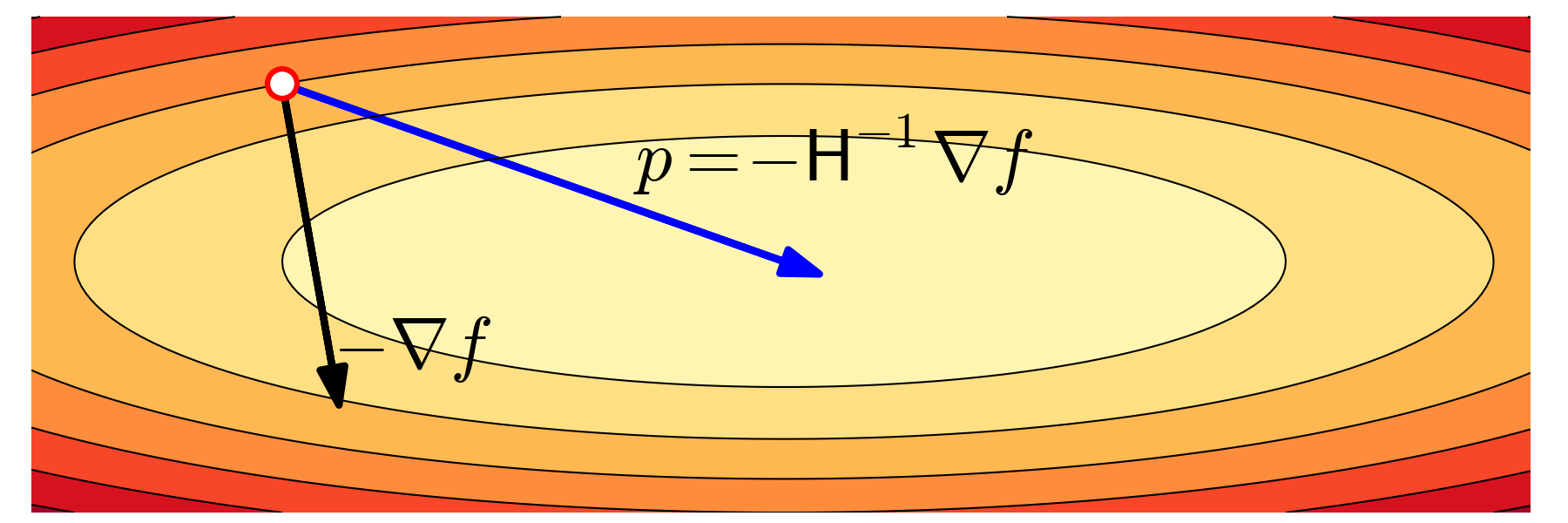}}
  \subfigure[inaccurate gradient]{%
        \includegraphics[width=0.49\textwidth]{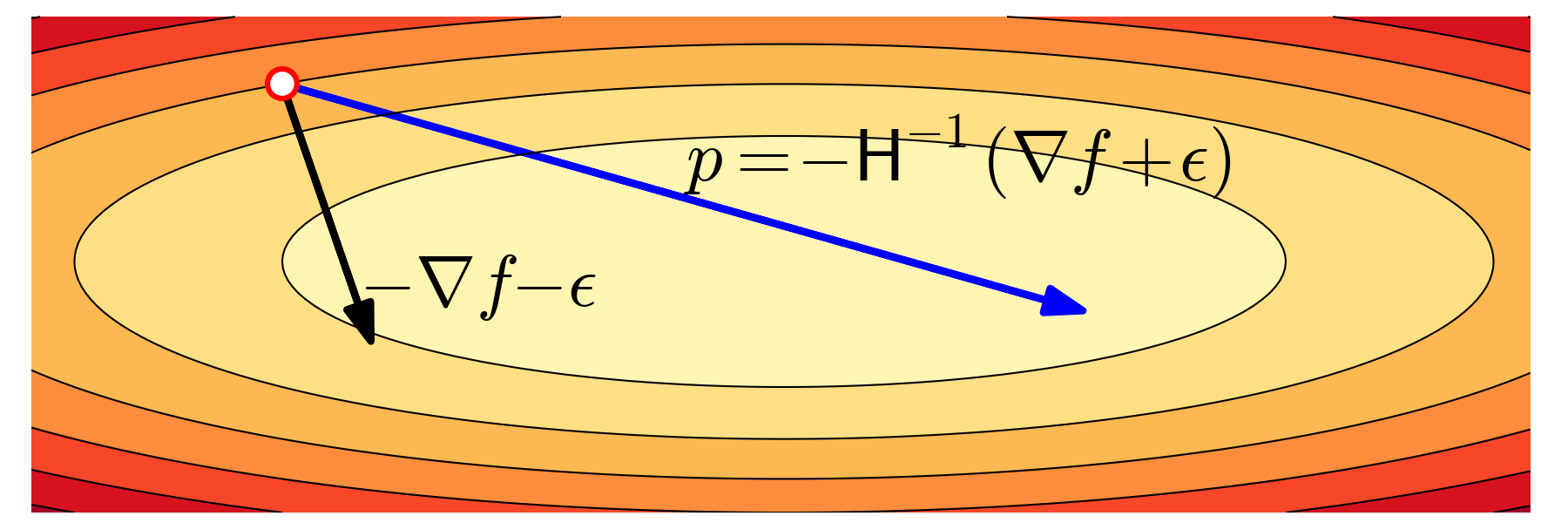}}
  \end{center}
  \caption{Even when the Hessian is exact, a small error in the gradient can
    significantly impact the quality of the step.  \label{fig:gradient_accuracy}}
\end{figure}

A na\"ive choice would be to set $\bar{g} = \nabla f(x_{0})$, \ie the gradient
evaluated at the initial point\footnote{The initial point for Arnoldi sampling
  is the current iterate of the outer optimization loop.}.  Although this is a
natural choice for accurate data, it can lead to inaccurate steps in the
presence of errors, even if the Hessian model is exact; this is illustrated in
Figure~\ref{fig:gradient_accuracy}.  In the following subsections, we propose
two possible choices for $\bar{g}$ that attempt to ameliorate errors in the
gradient.

\subsection{Variant 1: Step Average}

We can view each of the points $x_{j}, j=0,1,\ldots,m$, produced by Arnoldi
sampling as a potential location at which to center the quadratic model.  For
each of these locations, we can adopt the corresponding gradient, $g_{j}$, in
the model $q(p)$ and find a step.  Let $p_{j} = \mat{V} y_{j}$ denote this step.
Substituting this step into \eqref{eq:model}, and ignoring the trust-radius
constraint for the moment, the solution is
\begin{equation*}
  p_{j} = -\mat{V}\Lambda^{-1} \mat{V}^{T} g_{j}.
\end{equation*}
The trial solution for point $j$ is $x_{j} + p_{j}$.  Averaging over all potential
trial solutions, we obtain
\begin{equation*}
  \xnew = \frac{1}{m+1} \sum_{j=0}^{m} (x_{j} + p_{j}) 
  = \bar{x} - \mat{V}\Lambda^{-1} \mat{V}^{T} \bar{g},
\end{equation*}
where $\bar{x} = (\sum_{j} x_{j})/(m+1)$ and $\bar{g} = (\sum_{j} g_{j})/(m+1)$.
In the sequel, we will refer to this as the step-average approach.  Note that,
if the trust-radius constraint is present, we solve~\eqref{eq:model} with
$\bar{g}$ defined as the average gradient.

\subsection{Variant 2: Directional Derivatives}

Our numerical experiments suggest that step averaging is effective when the
error has mean zero; however, it does not address bias in the error.  This is
not surprising, because $\bar{g}$ is a simple average, so the bias will persist.
As an aside, we remark that constant bias in the set $\{g_{j}\}_{j=0}^{m}$ does
not impact the eigenvalue and eigenvector approximation in Arnoldi sampling,
because these approximations are based on differences in the gradient.

To address bias in the construction of $\bar{g}$, we examine the form of the
(unconstrained) solution when the step is in the subspace $\mat{V}$ and
$\bar{g}$ is the exact gradient:
\begin{equation*}
  p = \mat{V} y,
  \qquad\text{where}\qquad 
  y = \Lambda^{-1} \left(-\mat{V}^{T} g \right).
\end{equation*}
We see that the reduced-space solution, $y$, is the inverse of $\Lambda$ acting
on the reduced-space negative gradient.  This is an important observation,
because it means we can focus on approximating $\mat{V}^{T} g$ rather than $g$.

The reduced gradient, $\mat{V}^{T} g$, is the directional derivative of the
objective in the estimated eigenvector directions, $\mat{V}$.  These directional
derivatives can be approximating using the sampled function values:
\begin{equation*}
  \mat{V}^{T} g =
  \tilde{\mat{V}}_{1:m,1:r}^{T} \mat{Z}_{m}^{T} g \approx 
  \gred \equiv 
  \frac{1}{\alpha} \tilde{\mat{V}}_{1:m,1:r}^{T}
  \begin{bmatrix} 
    f_{1} - f_{0} \\ f_{2} - f_{0} \\ \vdots \\ f_{r} - f_{0}
  \end{bmatrix},
\end{equation*}
where $\gred$ denotes the approximate reduced gradient.  By using the function
values to approximate the directional derivatives, we eliminate bias in the
gradient error.

When the trust-region constraint is present, the reduced-space solution is found
by solving the trust-region problem \eqref{eq:model} with $\bar{g} = \mat{V}
\gred$, and adding the step, $p= \mat{V}y$, to the point $x_{0}$ as defined in
Arnoldi sampling.  We will refer to this as the directional-derivative variant.

\subsection{Comparison of the Two Variants}

We use the synthetic quadratic objective, defined earlier
by~\eqref{eq:quad_obj}, to study and compare the two proposed variants for
$\bar{g}$.  As before, we use $n = 256$ variables and set $\alpha = 1$.  The
initial point provided to the Arnoldi sampling algorithm is $(x_{0})_{i} =
\sin(i), i = 1,\ldots,256$.  The Hessian model is constructed using the $r = 4$
largest estimated eigenvalues and eigenvectors from Arnoldi sampling, which is
run with $m=16$ iterations.

\begin{figure}[tb]
  \begin{center}
  \subfigure[results when gradient error has no bias \label{fig:no_bias}]{%
    \includegraphics[width=0.9\textwidth]{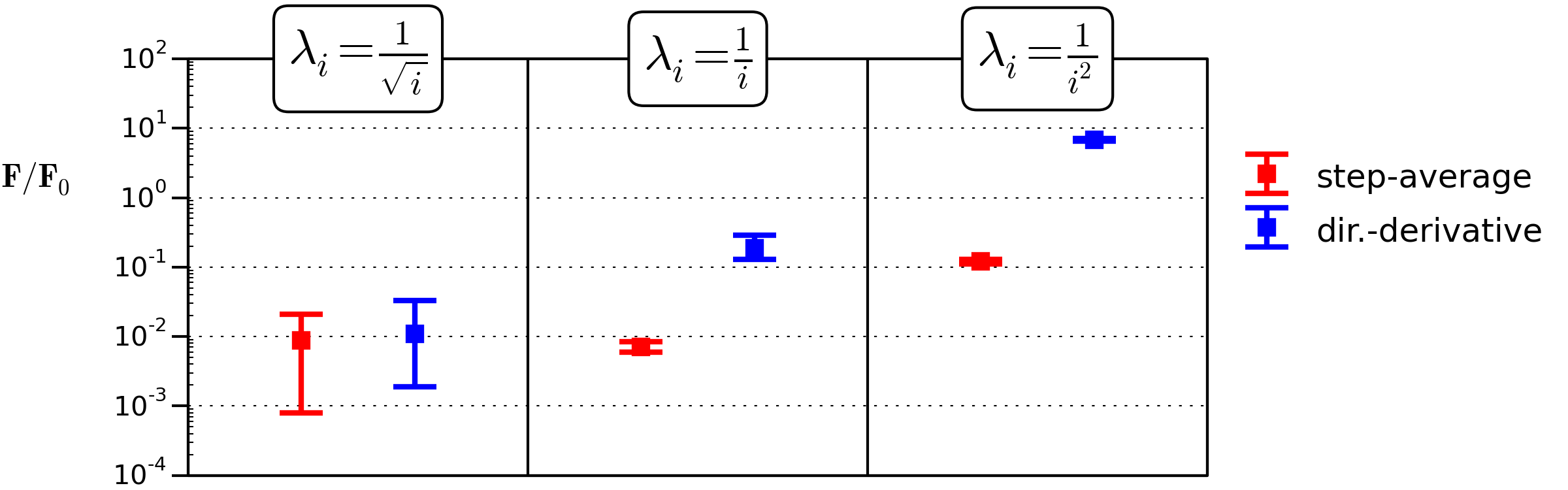}}\\
  \subfigure[results when gradient error has bias \label{fig:with_bias}]{%
        \includegraphics[width=0.9\textwidth]{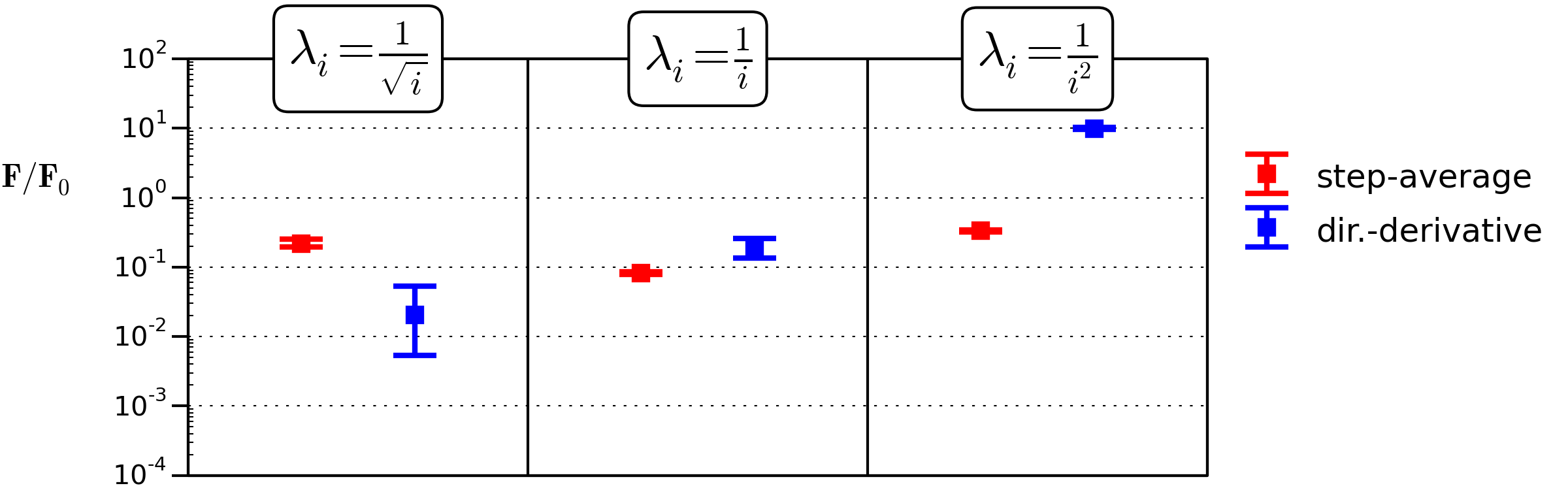}}
  \end{center}
  \caption{Comparison of the two methods of computing $\bar{g}$ --- the
    step-average and directional-derivative approaches --- without bias in the
    gradient error (figure a) and with bias in the gradient error (figure b).
  \label{fig:variants}}
\end{figure}

We consider two models for error.  The first model adds Gaussian noise with mean
zero to the function, $F$, and its gradient, $\nabla F$.  The noise has a
standard deviation equal to 2.5\% of $F(x_{0})$ and $\| \nabla F(x_{0}) \|$ for
the function and gradient components, respectively.  The second model for the
error is also Gaussian and has the same standard deviation, but each component
of the gradient error has a mean of $0.1 \| \nabla F(x_{0}) \| $, \ie the
gradient has a biased error.

The two variants of $\bar{g}$ were used to solve the
subproblem~\eqref{eq:model}, based on $F(x)$, for the three eigenvalue spectra,
$\{1/\sqrt{i}\}_{i=1}^{256}$, $\{1/i\}_{i=1}^{256}$, and
$\{1/i^2\}_{i=1}^{256}$.  For each spectra, the methods were applied 1000 times,
and statistics for the relative change in the objective, $F(x_0 + p)/F(x_{0})$,
were gathered.  Figure~\ref{fig:no_bias} shows the median (square symbol) for
the relative change in the objective when the error has no bias.  The bars show
the 0.025 and 0.975 quantiles of probability.  Figure~\ref{fig:with_bias} plots
the same results when bias is present.

When the error has no bias, it is clear that the step-average approach
outperforms the directional-derivative approach for this set of quadratics.
Indeed, for spectra with rapid decay, the objective increases on average using
the directional-derivative approach.  It must be emphasized, again, that these
results are sensitive to the choice of $\alpha$; the directional-derivative
variant relies on a finite-difference step-size that is accurate for both the
Hessian-vector product and gradient approximations.  When the objective function
varies rapidly in some directions, as it does for $\lambda_{i} = 1/i^2$, a good
step size for both the Hessian-vector products and gradient is difficult to find
and, in fact, may not exist.

As expected, the results change when the gradient error has bias.  In
particular, we see that the directional-derivative approach is relatively
insensitive to the addition of bias, \ie for a given spectra, the results are
similar.  In contrast, the step-average approach performs more poorly when the
spectrum decays slowly; however, the step-average approach continues to
outperform the directional-derivative variant when the spectrum decays rapidly.

\subsection{Optimization Framework}

When the objective function is nonlinear, the optimization algorithm must use
some form of globalization.  To this end, we have adopted a trust-region
framework~\cite{conn:2000} for this work.  We will refer to the overall
optimization algorithm, listed in Algorithm~\ref{alg:SAM}, as the Stochastic
Arnoldi's Method, or SAM for short.  Note that there are two variants of SAM; one
based on the step-average and one based on the directional-derivative method to
compute $\bar{g}$.

Each iteration of SAM begins by assessing convergence using the quadratic
model's gradient norm, $\| \bar{g} \|$.  In the case of the step-average
variant, the criterion is based on the average gradient norm.  The
directional-derivative variant uses the gradient norm of $\gred$.

If the convergence criterion is not met, SAM solves the (small) trust-region
subproblem~\eqref{eq:model} using the Mor\'e and Sorensen
algorithm~\cite{more:1983}.  In Algorithm~\ref{alg:SAM}, this subproblem is
written in terms of the reduced-space $y$; note that the step length satisfies
$\|p\| = \| \mat{V} y \| = \| y\|$, because the approximate eigenvectors are
orthonormal.

Subsequently, the step is used to obtain the trial solution $\xnew$.  The
trial-solution formula depends on whether the step-average or
directional-derivative variant is adopted.  Once the trial solution has been
computed, SAM follows a fairly standard trust-radius update based on the ratio,
$\rho$, of the actual objective reduction to the model's prediction of the
reduction.  One departure from conventional trust-region methods arises if the
step is rejected: in this case the proposed method re-evaluates the objective
and gradient at the current solution estimate.

While trust-region algorithms are well studied for accurate and even
inexact\footnote{Inexact refers to errors that can be controlled, \eg errors due
  to incomplete convergence of a solver.} function evaluations, their use with
imperfect data is less developed.  We plan a careful study of the theoretical
properties of Algorithm~\ref{alg:SAM} in future work, but make no guarantees
regarding its convergence at present.

\begin{algorithm}[tbp]\DontPrintSemicolon
  \KwData{$x_{0}$, $r$ (approximate-Hessian rank), $m$ (Arnoldi iterations),
    $\alpha$ (sample radius), $\Delta$ (initial trust radius), and $\tau$
    (convergence tolerance)}
  \KwResult{approximate minimum, $x$}
  \BlankLine
  set $x \leftarrow x_{0}$, and compute $\mat{F}_{1,1} = f(x)$ and 
  $\mat{G}_{:,1} = g(x)$\;
  use Arnoldi sampling (Algorithm~\ref{alg:Arnoldi_sample}) to obtain 
  $\mat{X}$, $\mat{F}$, $\mat{G}$, $\mat{V}$, and $\Lambda$\;
  \For{$k = 1,2,\ldots$}{
    \If(\algcmnt{check convergence}){%
      $\| \bar{g} \| \leq \tau$}{
      return\;
    } 
    Solve trust-region optimization problem in the reduced space:
    \begin{equation*}\setlength{\abovedisplayskip}{1ex}\setlength{\belowdisplayskip}{-1ex}
      \min_{y} \; q(y) = \bar{f} + \bar{g}^{T} \mat{V} y + \frac{1}{2} y^{T} \Lambda y, \qquad \text{subject to}\; \| y \| \leq \Delta.
    \end{equation*}\;
    compute $\xnew \leftarrow \bar{x} + \mat{V}y$ (step-average) 
    or $\xnew \leftarrow x + \mat{V} y$ (directional-derivative)\;
    compute $f_{k} = f(\xnew)$ and $g_{k} = g(\xnew)$\;
    compute $\rho = \textsf{ared}/\textsf{pred} = -(f_{k-1} - f_{k})/q(y)$\;
    \eIf{$\rho < 0.1$}{
      $\Delta \leftarrow \Delta/4$\;
    }{%
      \If{$\rho > 3/4$ and $\|p\| < \Delta$}{
        $\Delta \leftarrow \min(2\Delta, \Delta_{\max})$\;
      }
    }
    \eIf{$\rho > 10^{-4}$}{%
      set $x \leftarrow \xnew$, and update $\mat{F}_{:,1} = f_{k}$ and $\mat{G}_{:,1} = g_{k}$\;
    }{
      keep $x$ and resample $\mat{F}_{:,1} = f(x)$ 
      and $\mat{G}_{:,1} = g(x)$\;      
    }
    use Arnoldi sampling (Algorithm~\ref{alg:Arnoldi_sample}) to update $\mat{X}_{:,2:(m+1)}$, 
    $\mat{F}_{1,2:(m+1)}$, $\mat{G}_{:,2:(m+1)}$, and $\mat{V}_{:,1:r}$\;
  }
  \caption{Stochastic Arnoldi's Method.\label{alg:SAM}}
\end{algorithm}

\section{Results}

We conclude by benchmarking SAM against two established optimization algorithms:
the (derivative-based) BFGS~\cite{nocedal:2006} quasi-Newton method and the
(derivative-free) Nelder-Mead algorithm~\cite{nelder:1965}.  For this study, we
optimize a modified multi-dimensional Rosenbrock function:
\begin{equation*}
  F(x) = \sum_{i=1}^{n/2} \frac{1}{i}\left[100 (x_{2i} - x_{2i-1}^{2})^2 + (1 -
    x_{2i-1})^2 \right],
\end{equation*}
where $n = 256$.  This definition differs from the standard multi-dimensional
Rosenbrock function~\cite{dixon:1994} with the introduction of the scaling factor
$1/i$.  This factor ensures the Hessian has a decaying spectrum, which is an
underlying assumption in Arnoldi sampling and, consequently, SAM.

As before, we consider both unbiased and biased Gaussian noise with standard
deviations of 2.5\% of $F(x_{0})$ and $\| \nabla F(x_{0}) \|$ for the function
and gradient components, respectively.  For the biased-noise case, each entry in
the gradient has a mean error of $0.1 \| \nabla F(x_{0}) \|$.

We benchmark the two variants of SAM against the BFGS and Nelder-Mead
implementations in Matlab\textsuperscript{\textregistered}.  The initial iterate
for SAM and BFGS is given by
\begin{equation*}
  \left(x_{0}\right)_{i} = \begin{cases}
    -1, &\text{if $i$ is odd}, \\
    0, &\text{if $i$ is even}.
  \end{cases}
\end{equation*}
For Nelder-Mead, $x_{0}$ forms one of the vertices of the initial simplex.  The
parameters used for SAM are as follows: $r = 4$, $m=16$, $\alpha =0.5$, $\Delta
= 10 \|x_{0}\|$, and $\tau = 0.1$.  In addition, SAM is limited to a maximum of
10 iterations.

\begin{figure}[tbp]
  \begin{center}
    \includegraphics[width=\textwidth]{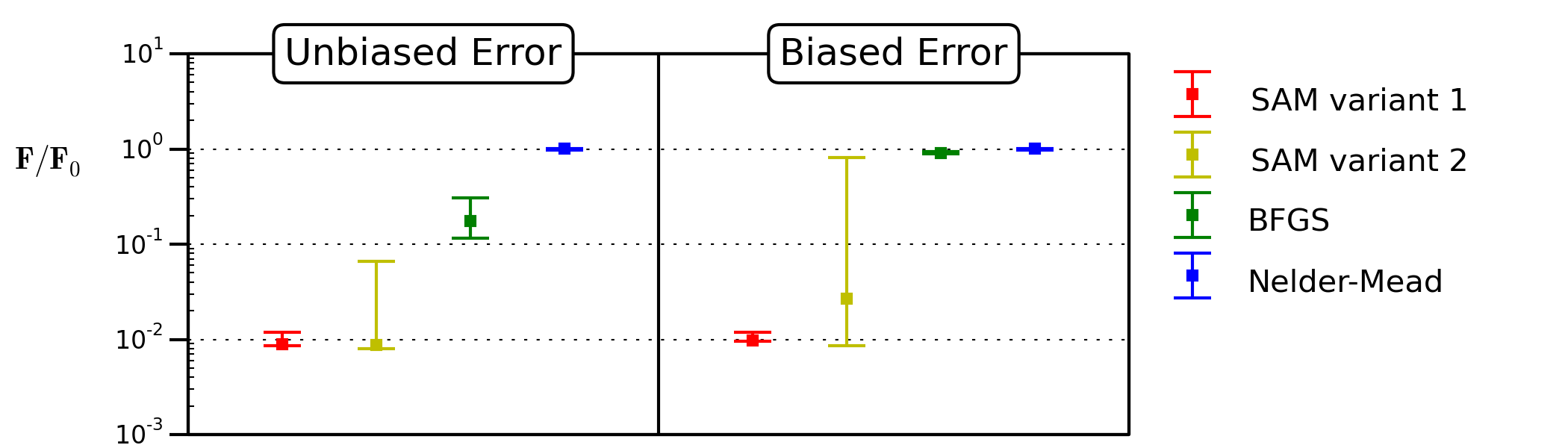}
    \caption{Comparison of the SAM variants with BFGS and Nelder-Mead on the
      multi-dimensional Rosenbrock function with Gaussian noise added.
      \label{fig:rosenbrock_compare}}
  \end{center}
\end{figure}

Figure~\ref{fig:rosenbrock_compare} compares the median objective reduction
achieved by the two SAM variants, BFGS, and Nelder-Mead.  The lower and upper
bars denote the 0.025 and 0.975 probability quantiles, respectively.  The
step-average variant of SAM reduces the objective by approximately two orders of
magnitude for both the unbiased and biased errors.  Somewhat surprisingly, the
step-average method performs better than the directional-derivative method in
the biased-error case; this may be due to the choice of model radius $\alpha$.
Both variants of SAM perform significantly better than the BFGS and Nelder-Mead
simplex methods.

\begin{figure}[tbp]
  \begin{center}
    \subfigure[Unbiased noise]{%
      \includegraphics[width=0.47\textwidth]{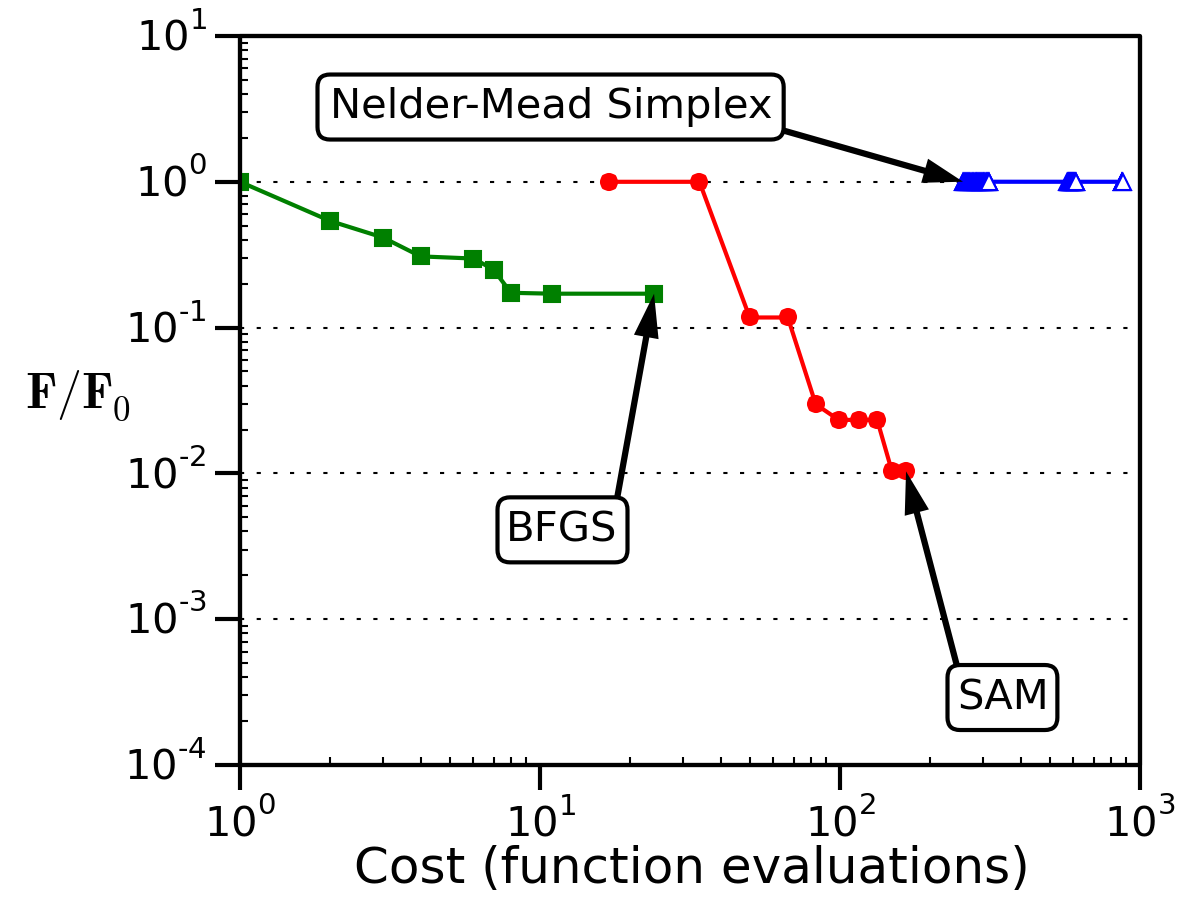}}\hfill
    \subfigure[Biased noise]{%
      \includegraphics[width=0.47\textwidth]{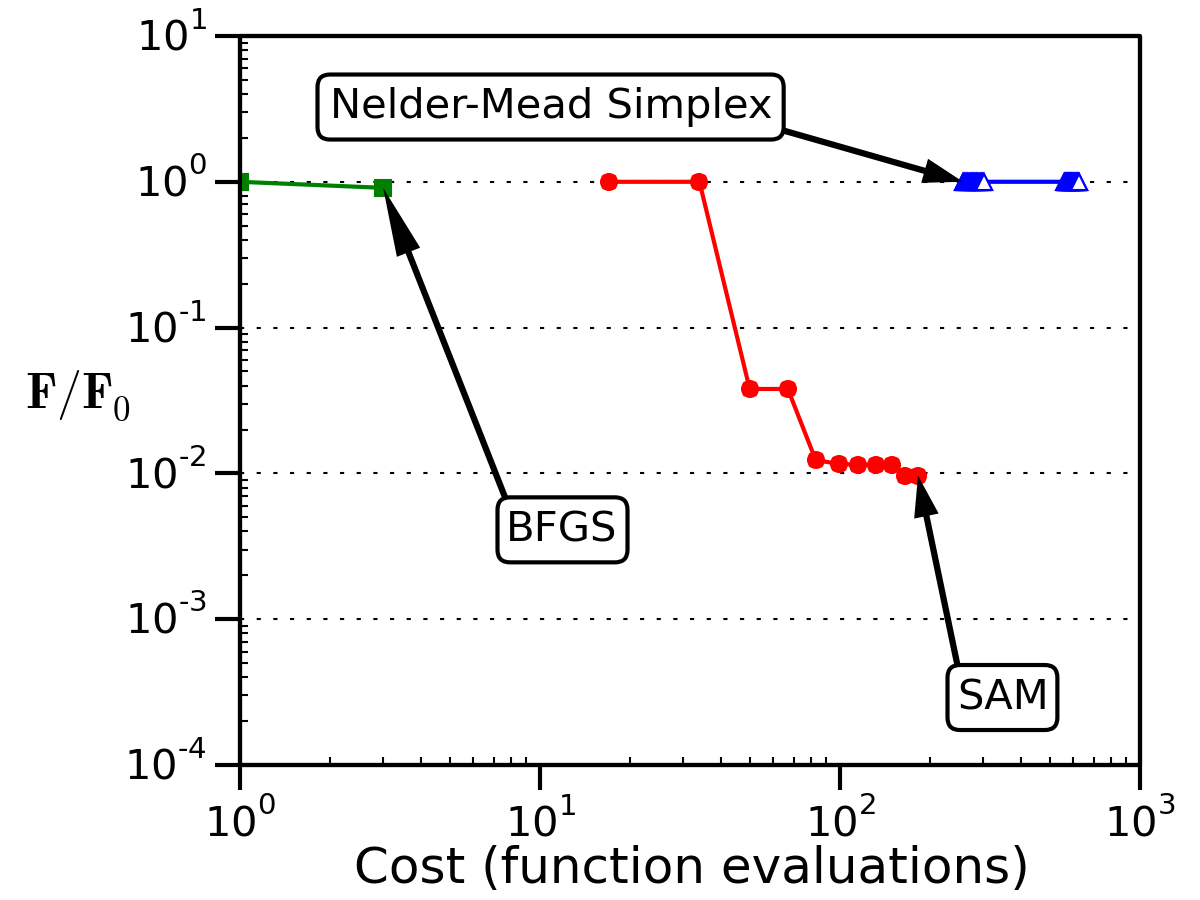}}
    \caption{Sample convergence histories for the SAM (step-average variant),
      BFGS, and Nelder-Mead algorithms.\label{fig:converge}}
  \end{center}
\end{figure}

Figure~\ref{fig:converge} compares the convergence histories of SAM with those
of BFGS and Nelder-Mead.  Only the step-average variant of SAM is shown.  The
BFGS convergence stalls in both cases, because the errors cause the line search
to fail.  These histories, which are typical, demonstrate that the performance
of SAM comes at the cost of additional functional evaluations.

\section{Summary and Conclusions}

Numerical optimization has helped engineers in numerous fields; however, there
remain important applications that cannot use conventional optimization
algorithms.  In this work, we have targeted applications that have
large-dimensional design spaces and whose outputs are imperfect, \ie their
outputs and derivatives contain irreducible errors that are incompatible with
most gradient-based algorithms.  To enable optimization for these applications, we
have developed a high-dimensional sampling method based on Arnoldi's method.

Arnoldi sampling adaptively selects new sample points and tends to capture the
dominant eigenmodes of the (spatially averaged) Hessian.  The sampling strategy
could be used in conjunction with a nonparametric regression, \eg
Gaussian-process regression, to build surrogate models for optimization when
errors are present in the objective and its derivatives; however, in the present
work, we have focused on conventional quadratic models.

We investigated two methods of constructing the linear term that appears in the
quadratic model.  The first approach involved averaging over all possible sample
points and gradients.  The second approach used directional derivatives of the
objective function to estimate the gradient in the reduced space.  The
directional-derivative variant was found to be effective when biased errors were
present and the Hessian's spectrum was slowly decaying; however, in general, the
step-average approach was more effective.

A trust-region algorithm called the Stochastic Arnoldi's Method (SAM) was
proposed and used to minimize the multi-dimensional Rosenbrock function.
Synthetic errors were introduced in the evaluation of the objective and its
gradient.  The performance of SAM on this problem was promising relative to
derivative-based (BFGS) and derivative-free (Nelder-Mead) algorithms.  In
particular, whereas BFGS had inconsistent convergence and Nelder-Mead failed to
converge, the Arnoldi-based algorithm consistently reduced the objective by two
orders of magnitude.

\bibliographystyle{aiaa}
\bibliography{/home/jehicken/Biblio/jehicken}

\begin{thebibliography}{10}
\newcommand{\enquote}[1]{``#1''}

\bibitem{atag:2014}
\enquote{{Air Transport Action Group: Facts \& Figures},}
  http://www.atag.org/facts-and-figures.html, 2014, accessed: 07/14/2014.

\bibitem{spalart:2009}
Spalart, P.~R., \enquote{{Detached-eddy simulation},} {\em Annual Review of
  Fluid Mechanics\/}, Vol.~41, 2009, pp.~181--202.

\bibitem{debaar:2005}
De~Baar, M.~R., Thyagaraja, A., Hogeweij, G. M.~D., Knight, P.~J., and Min, E.,
  \enquote{{Global plasma turbulence simulations of q=3 sawtoothlike events in
  the RTP tokamak},} {\em Physical review letters\/}, Vol.~94, No.~3, 2005,
  pp.~035002+.

\bibitem{enotiadis:1990}
Enotiadis, A.~C., Vafidis, C., and Whitelaw, J.~H., \enquote{{Interpretation of
  cyclic flow variations in motored internal combustion engines},} {\em
  Experiments in fluids\/}, Vol.~10, No. 2-3, 1990, pp.~77--86.

\bibitem{lorenz:1963}
Lorenz, E.~N., \enquote{{Deterministic nonperiodic flow},} {\em Journal of the
  Atmospheric Sciences\/}, Vol.~20, No.~2, 1963, pp.~130--141.

\bibitem{lea:2000}
Lea, D.~J., Allen, M.~R., and Haine, T. W.~N., \enquote{{Sensitivity analysis
  of the climate of a chaotic system},} {\em Tellus A\/}, Vol.~52, No.~5, Oct.
  2000, pp.~523--532.

\bibitem{eyink:2004}
Eyink, G.~L., Haine, T. W.~N., and Lea, D.~J., \enquote{{Ruelle's linear
  response formula, ensemble adjoint schemes and L\'{e}vy flights},} {\em
  Nonlinearity\/}, Vol.~17, No.~5, Sept. 2004, pp.~1867+.

\bibitem{ashley:2014}
Ashley, A. and Hicken, J.~E., \enquote{{Optimization Algorithm for Systems
  Governed by Chaotic Dynamics},} {\em 2014 AIAA Aviation Conference\/}, June
  2014, {AIAA} 2014-2434.

\bibitem{wang:2013b}
Wang, Q., Hu, R., and Blonigan, P., \enquote{{Sensitivity computation of
  periodic and chaotic limit cycle oscillations},} Aug. 2013,
  arXiv:1204.0159v4.

\bibitem{veroy:2005}
Veroy, K. and Patera, A.~T., \enquote{{Certified real-time solution of the
  parametrized steady incompressible Navier–Stokes equations: rigorous
  reduced-basis a posteriori error bounds},} {\em Int. J. Numer. Meth.
  Fluids\/}, Vol.~47, No. 8-9, March 2005, pp.~773--788.

\bibitem{buithanh:2008}
Bui-Thanh, T., Willcox, K., and Ghattas, O., \enquote{{Model Reduction for
  Large-Scale Systems with High-Dimensional Parametric Input Space},} {\em SIAM
  Journal on Scientific Computing\/}, Vol.~30, No.~6, Jan. 2008,
  pp.~3270--3288.

\bibitem{bashir:2008}
Bashir, O., Willcox, K., Ghattas, O., van Bloemen~Waanders, B., and Hill, J.,
  \enquote{{Hessian-based model reduction for large-scale systems with
  initial-condition inputs},} {\em Int. J. Numer. Meth. Engng.\/}, Vol.~73,
  No.~6, Feb. 2008, pp.~844--868.

\bibitem{lockwood:2012}
Lockwood, B., Anitescu, M., and Mavriplis, D.~J., \enquote{{Mixed
  aleatory/epistemic uncertainty quantification for hypersonic flows via
  gradient-based optimization and surrogate models},} {\em 50th AIAA Aerospace
  Sciences Meeting\/}, Nashville, Tennessee, 2012, pp. 9--12.

\bibitem{boopathy:2014}
Boopathy, K. and Rumpfkeil, M.~P., \enquote{{Robust Optimizations of Structural
  and Aerodynamic Designs},} {\em 2014 AIAA Aviation Conference\/}, June 2014.

\bibitem{aftosmis:1997}
Aftosmis, M.~J., \enquote{Lecture notes for the 28th computational fluid
  dynamics lecture series: solution adaptive {C}artesian grid methods for
  aerodynamic flows with complex geometries,} Tech. rep., von {K}\'{a}rm\'{a}n
  {I}nstitute for {F}luid {D}ynamics, Rhode-Saint-Gen\`{e}se, Belgium, March
  1997.

\bibitem{aftosmis:1998}
Aftosmis, M.~J., Berger, M.~J., and Melton, J.~E., \enquote{{Robust and
  efficient Cartesian mesh generation for component-based geometry},} {\em AIAA
  journal\/}, Vol.~36, No.~6, 1998, pp.~952--960.

\bibitem{fidkowski:2007}
Fidkowski, K.~J. and Darmofal, D.~L., \enquote{A triangular cut-cell adaptive
  method for high-order discretizations of the compressible {N}avier-{S}tokes
  equations,} {\em Journal of Computational Physics\/}, Vol.~225, Aug. 2007,
  pp.~1653--1672.

\bibitem{peskin:1972}
Peskin, C.~S., \enquote{{Flow Patterns Around Heart Values: A Numerical
  Method},} {\em Journal of Computational Physics\/}, Vol.~10, No.~2, 1972,
  pp.~252--271.

\bibitem{mittal:2005}
Mittal, R. and Iaccarino, G., \enquote{{Immersed boundary methods},} {\em
  Annual Review of Fluid Mechanics\/}, Vol.~37, 2005, pp.~239--261.

\bibitem{griffith:2005}
Griffith, B.~E. and Peskin, C.~S., \enquote{{On the order of accuracy of the
  immersed boundary method: higher order convergence rates for sufficiently
  smooth problems},} {\em Journal of Computational Physics\/}, Vol.~208, No.~1,
  2005, pp.~75--105.

\bibitem{gunzburger:2003}
Gunzburger, M.~D., {\em {Perspectives in flow control and optimization}\/},
  Society for Industrial and Applied Mathematics, 2003.

\bibitem{giles:2000}
Giles, M.~B. and Pierce, N.~A., \enquote{{An introduction to the adjoint
  approach to design},} {\em Flow, Turbulence and Combustion\/}, Vol.~65,
  No.~3, 2000, pp.~393--415.

\bibitem{collis:2002}
Collis, S.~S. and Heinkenschloss, M., \enquote{Analysis of the streamline
  upwind/{P}etrov {G}alerkin method applied to the solution of optimal control
  problems,} Tech. Rep. TR02-01, Houston, Texas, 2002.

\bibitem{hicken:dual2014}
Hicken, J.~E. and Zingg, D.~W., \enquote{{Dual consistency and functional
  accuracy: a finite-difference perspective},} {\em Journal of Computational
  Physics\/}, Vol.~256, Jan. 2014, pp.~161--182.

\bibitem{mohammadi:2004}
Mohammadi, B. and Pironneau, O., \enquote{{Shape Optimization in Fluid
  Mechanics},} {\em Annual Review of Fluid Mechanics\/}, Vol.~36, No.~1, 2004,
  pp.~255--279.

\bibitem{derakhshan:2008}
Derakhshan, S., Mohammadi, B., and Nourbakhsh, A., \enquote{{Incomplete
  sensitivities for 3D radial turbomachinery blade optimization},} {\em
  Computers \& Fluids\/}, Vol.~37, No.~10, 2008, pp.~1354--1363.

\bibitem{balabanov:2004}
Balabanov, V. and Venter, G., \enquote{{Multi-fidelity optimization with
  high-fidelity analysis and low-fidelity gradients},} {\em 10th AIAA/ISSMO
  Multidisciplinary Analysis and Optimization Conference\/}, Albany, New York,
  2004.

\bibitem{conn:1998}
Conn, A.~R., Scheinberg, K., and Toint, P.~L., \enquote{{A derivative free
  optimization algorithm in practice},} {\em 7th AIAA/USAF/NASA/ISSMO Symposium
  on Multidisciplinary Analysis and Optimization\/}, St. Louis, Missouri, 1998.

\bibitem{powell:2002}
Powell, M. J.~D., \enquote{{UOBYQA: unconstrained optimization by quadratic
  approximation},} {\em Mathematical Programming\/}, Vol.~92, No.~3, 2002,
  pp.~555--582.

\bibitem{conn:2009}
Conn, A.~R., Scheinberg, K., and Vicente, L.~N., {\em {Introduction to
  Derivative-Free Optimization}\/}, Society for Industrial and Applied
  Mathematics, Jan. 2009.

\bibitem{nelder:1965}
Nelder, J.~A. and Mead, R., \enquote{{A simplex method for function
  minimization},} {\em The Computer Journal\/}, Vol.~7, No.~4, Jan. 1965,
  pp.~308--313.

\bibitem{holland:1975}
Holland, J.~H., {\em {Adaptation in Natural and Artificial Systems}\/}, The
  University of Michigan Press, Ann Arbor, Michigan, 1975.

\bibitem{hajela:1990}
Hajela, P., \enquote{{Genetic search --- an approach to the nonconvex
  optimization problem},} {\em AIAA Journal\/}, Vol.~28, No.~7, July 1990,
  pp.~1205--1210.

\bibitem{griewank:2008}
Griewank, A. and Walther, A., {\em {Evaluating derivatives: principles and
  techniques of algorithmic differentiation}\/}, Society for Industrial and
  Applied Mathematics, 2008.

\bibitem{robbins:1951}
Robbins, H. and Monro, S., \enquote{{A stochastic approximation method},} {\em
  The Annals of Mathematical Statistics\/}, 1951, pp.~400--407.

\bibitem{kiefer:1952}
Kiefer, J. and Wolfowitz, J., \enquote{{Stochastic estimation of the maximum of
  a regression function},} {\em The Annals of Mathematical Statistics\/},
  Vol.~23, No.~3, 1952, pp.~462--466.

\bibitem{spall:1992}
Spall, J.~C., \enquote{{Multivariate stochastic approximation using a
  simultaneous perturbation gradient approximation},} {\em IEEE Transactions on
  Automatic Control\/}, Vol.~37, No.~3, 1992, pp.~332--341.

\bibitem{spall:2005}
Spall, J.~C., {\em {Introduction to stochastic search and optimization:
  estimation, simulation, and control}\/}, John Wiley \& Sons, Hoboken, New
  Jersey, 2003.

\bibitem{spall:2009}
Spall, J.~C., \enquote{{Feedback and weighting mechanisms for improving
  Jacobian estimates in the adaptive simultaneous perturbation algorithm},}
  {\em IEEE Transactions on Automatic Control\/}, Vol.~54, No.~6, 2009,
  pp.~1216--1229.

\bibitem{pearson:1901}
Pearson, K., \enquote{{On lines and planes of closest fit to systems of points
  in space},} {\em The London, Edinburgh, and Dublin Philosophical Magazine and
  Journal of Science\/}, Vol.~2, No.~11, 1901, pp.~559--572.

\bibitem{willcox:2002}
Willcox, K. and Peraire, J., \enquote{{Balanced model reduction via the proper
  orthogonal decomposition},} {\em AIAA journal\/}, Vol.~40, No.~11, 2002,
  pp.~2323--2330.

\bibitem{rowley:2004}
Rowley, C.~W., Colonius, T., and Murray, R.~M., \enquote{{Model reduction for
  compressible flows using POD and Galerkin projection},} {\em Physica D:
  Nonlinear Phenomena\/}, Vol.~189, No.~1, 2004, pp.~115--129.

\bibitem{volkwein:2011}
Volkwein, S., \enquote{{Model reduction using proper orthogonal
  decomposition},} Lecture Notes, Institute of Mathematics and Scientific
  Computing, University of Graz, 2011,
  http://www.math.uni-konstanz.de/numerik/personen/volkwein/teaching/POD-Vorlesung.pdf.

\bibitem{amsallem:2012}
Amsallem, D. and Farhat, C., \enquote{{Stabilization of projection-based
  reduced-order models},} {\em International Journal for Numerical Methods in
  Engineering\/}, Vol.~91, No.~4, 2012, pp.~358--377.

\bibitem{lieberman:2010}
Lieberman, C., Willcox, K., and Ghattas, O., \enquote{{Parameter and state
  model reduction for large-scale statistical inverse problems},} {\em SIAM
  Journal on Scientific Computing\/}, Vol.~32, No.~5, 2010, pp.~2523--2542.

\bibitem{mckay:1979}
McKay, M.~D., Beckman, R.~J., and Conover, W.~J., \enquote{{Comparison of three
  methods for selecting values of input variables in the analysis of output
  from a computer code},} {\em Technometrics\/}, Vol.~21, No.~2, 1979,
  pp.~239--245.

\bibitem{saad:1993}
Saad, Y., \enquote{A flexible inner-outer preconditioned {GMRES} algorithm,}
  {\em {SIAM} Journal on Scientific and Statistical Computing\/}, Vol.~14,
  No.~2, 1993, pp.~461--469.

\bibitem{saad:2003}
Saad, Y., {\em {Iterative Methods for Sparse Linear Systems}\/}, SIAM,
  Philadelphia, PA, 2nd ed., 2003.

\bibitem{wang:1992}
Wang, Z., Navon, I.~M., Dimet, F.~X., and Zou, X., \enquote{{The second order
  adjoint analysis: Theory and applications},} {\em Meteorology and Atmospheric
  Physics\/}, Vol.~50, 1992, pp.~3--20.

\bibitem{borzi:2011}
Borz\`{\i}, A. and Schulz, V., {\em {Computational Optimization of Systems
  Governed by Partial Differential Equations}\/}, Society for Industrial and
  Applied Mathematics, Jan. 2011.

\bibitem{hicken:inexact2014}
Hicken, J.~E., \enquote{{Inexact Hessian-vector products in reduced-space
  differential-equation constrained optimization},} {\em Optimization and
  Engineering\/}, Vol.~15, No.~3, Sept. 2014, pp.~575--608.

\bibitem{arnoldi:1951}
Arnoldi, W.~E., \enquote{{The principle of minimized iterations in the solution
  of the matrix eigenvalue problem},} {\em Quarterly of Applied Mathematics\/},
  Vol.~9, No.~1, 1951, pp.~17--29.

\bibitem{conn:2000}
Conn, A.~R., Gould, N. I.~M., and Toint, P.~L., {\em {Trust Region Methods}\/},
  Society for Industrial and Applied Mathematics, Jan. 2000.

\bibitem{more:1983}
Mor\'{e}, J.~J. and Sorensen, D.~C., \enquote{{Computing a Trust Region Step},}
  {\em SIAM Journal on Scientific and Statistical Computing\/}, Vol.~4, No.~3,
  Sept. 1983, pp.~553--572.

\bibitem{nocedal:2006}
Nocedal, J. and Wright, S.~J., {\em {Numerical Optimization}\/},
  Springer--Verlag, Berlin, Germany, 2nd ed., 2006.

\bibitem{dixon:1994}
Dixon, L. C.~W. and Mills, D.~J., \enquote{{Effect of rounding errors on the
  variable metric method},} {\em Journal of Optimization Theory and
  Applications\/}, Vol.~80, No.~1, 1994, pp.~175--179.

\end{thebibliography}

\end{document}